# CHAÎNES DE MARKOV CONSTRUCTIVES INDEXÉES PAR Z

By Jean Brossard and Christophe Leuridan

*Institut Fourier*

Nous étudions les chaînes de Markov $(X_n)_{n\in\mathbf{Z}}$ gouvernées par une relation de récurrence de la forme $X_{n+1} = f(X_n, V_{n+1})$, où $(V_n)_{n\in\mathbf{Z}}$ est une suite de variables aléatoires indépendantes et de même loi telle pour tout $n \in \mathbf{Z}$, $V_{n+1}$ est indépendante de la suite $((X_k, V_k))_{k\leq n}$. L'objet de l'article est de donner une condition nécessaire et suffisante pour que les innovations $(V_n)_{n\in\mathbf{Z}}$ déterminent complètement la suite $(X_n)_{n\in\mathbf{Z}}$ et de décrire l'information manquante dans le cas contraire.

**Introduction.** Dans cet article, nous étudions la filtration d'une chaîne de Markov constructive indexée par $\mathbf{Z}$. Nous appelons chaîne de Markov constructive (homogène) une suite $(X_n)_{n\in\mathbf{Z}}$ de variables aléatoires à valeurs dans un espace d'états $(E, \mathcal{E})$ satisfaisant une relation de récurrence de la forme $X_{n+1} = f(X_n, V_{n+1})$, où $(V_n)_{n\in\mathbf{Z}}$ est une suite de variables aléatoires indépendantes et de même loi à valeurs dans un espace $(G, \mathcal{G})$, $f$ est une application mesurable de $(E \times G, \mathcal{E} \otimes \mathcal{G})$ dans $(E, \mathcal{E})$ et $V_{n+1}$ est indépendante de la suite $((X_k, V_k))_{k\leq n}$ pour tout $n \in \mathbf{Z}$. Nous appelons innovations les variables aléatoires $V_n$ qui fournissent l'information nouvelle (indépendante du passé) à chaque instant. Sous ces conditions, la suite $(X_n)_{n\in\mathbf{Z}}$ est une chaîne de Markov dans la filtration naturelle de $((X_n, V_n))_{n\in\mathbf{Z}}$, de noyau de transition $\Pi$ où $\Pi(x, \cdot)$ est la loi de $f(x, V_1)$.

Les chaînes de Markov constructives fournissent beaucoup d'exemples de chaînes de Markov et apparaissent naturellement en simulation. La donnée d'une variable aléatoire $X_0$ et d'une suite $(V_n)_{n\geq 1}$ de variables i.i.d., indépendante de $X_0$ permet à construire la suite $(X_n)_{n\in\mathbf{Z}}$ vérifiant la relation de récurrence $X_{n+1} = f(X_n, V_{n+1})$. En revanche, pour les chaînes de Markov indexées par $\mathbf{Z}$, on ne dispose pas de condition initiale et en général la connaissance de $X_0$ et de la suite $(V_n)_{n\in\mathbf{Z}}$ ne permet pas de construire la suite $(X_n)_{n\in\mathbf{Z}}$.

Pour tout entier $N$ (aussi proche de $-\infty$ soit-il), la connaissance de $(X_n)_{n \leq N}$ et de la suite $(V_n)_{n \in \mathbf{Z}}$ détermine complètement la suite $(X_n)_{n \in \mathbf{Z}}$ par la relation de récurrence $X_{n+1} = f(X_n, V_{n+1})$. En notant $(\mathcal{F}_n^Y)_{n \in \mathbf{Z}}$ la filtration naturelle (complétée) d'une suite de variables aléatoires $(Y_n)_{n \in \mathbf{Z}}$, on a donc pour tout $n \geq N$, $\mathcal{F}_n^{(X,V)} = \mathcal{F}_N^X \vee \mathcal{F}_n^V$. Il est tentant de faire tendre $N$ vers $-\infty$ dans cette égalité et d'écrire $\mathcal{F}_n^{(X,V)} = \mathcal{F}_{-\infty}^X \vee \mathcal{F}_n^V$, comme l'ont fait Kallianpur et Wiener en 1956 dans un contexte semblable (voir [7] pour une discussion de ce point). En fait, il arrive que l'inclusion triviale

$$\left(\bigcap_{N \in \mathbf{Z}} \mathcal{F}_N^X\right) \vee \mathcal{F}_n^V \subset \bigcap_{N \in \mathbf{Z}} \mathcal{F}_N^X \vee \mathcal{F}_n^V$$

soit stricte et la propriété d'échange (de l'intersection décroissante sur $N$ avec le suprémum des tribus) n'a pas lieu en général. La propriété d'échange a été étudiée de façon générale par von Weizsäcker [12].

Un exemple très simple de ce phénomène, dû à Vinokourov mais non publié, est celui où $(X_n)_{n \in \mathbf{Z}}$ est une chaîne de Markov sur $\{-1, 1\}$ de probabilités de transition $p(1, 1) = p(-1, -1) = p$ et $p(-1, 1) = p(1, -1) = q$ avec $p$ et $q$ positifs de somme 1. Dans ce cas, on peut écrire $X_{n+1} = X_n V_{n+1}$ où $V_{n+1} = X_n X_{n+1}$ est une variable aléatoire indépendante de $\mathcal{F}_n^{(X,V)} = \mathcal{F}_n^X$ et de loi $p\delta_1 + q\delta_{-1}$. Bien que la tribu asymptotique $\mathcal{F}_{-\infty}^X$ soit triviale, l'inclusion triviale $\mathcal{F}_n^V \subset \mathcal{F}_n^{(X,V)}$ est stricte. En effet, par symétrie, la variable aléatoire $X_n$ est indépendante de $\mathcal{F}_n^V$ (et même de $\mathcal{F}_{+\infty}^V$), et de loi uniforme sur $\{-1, 1\}$. Et on vérifie immédiatement que la variable aléatoire $X_n$ fournit exactement l'information manquante, c'est-à-dire que $\mathcal{F}_n^{(X,V)} = \mathcal{F}_n^V \vee \sigma(X_n)$.

Au vu de cet exemple, deux problèmes se posent naturellement. Le premier est de déterminer à quelle condition la suite $V = (V_n)_{n \in \mathbf{Z}}$ détermine complètement la suite $X = (X_n)_{n \in \mathbf{Z}}$. Le second est de décrire l'information manquante dans le cas contraire.

Remarquons que si $V$ détermine $X$, alors la filtration $(\mathcal{F}_n^X)_{n \in \mathbf{Z}}$ est non seulement incluse, mais aussi immergée dans la filtration de type produit $(\mathcal{F}_n^V)_{n \in \mathbf{Z}}$: elle est donc standard (voir [3] pour la définition de ces notions). Pour que $V$ détermine $X$, il est donc nécessaire que la filtration $(\mathcal{F}_n^X)_{n \in \mathbf{Z}}$ soit standard, mais cela n'est aucunement suffisant: on vérifie facilement que dans l'exemple de Vinokourov, $(\mathcal{F}_n^X)_{n \in \mathbf{Z}}$ est standard.

L'objet de cet article n'est pas d'étudier les propriétés de standardité, mais plutôt les deux problèmes ci-dessus. L'article est organisé en trois parties.

Dans la première partie, nous nous plaçons dans le cadre d'une chaîne de Markov inhomogène gouvernée par une relation de récurrence $X_{n+1} = f_n(X_n, V_{n+1})$. Nous montrons que les tribus asymptotiques $\mathcal{F}_{-\infty}^{X,V}$ et $\mathcal{F}_{-\infty}^X$ sont égales et étudions de façon générale la loi conditionnelle de $X$ sachant $\sigma(V) \vee \mathcal{F}_{-\infty}^{X,V}$.

Dans la deuxième partie, nous nous plaçons dans le cadre d'une chaîne de Markov homogène sur un espace d'états dénombrable. Lorsque la chaîne $(X_n)_{n\in\mathbf{Z}}$ est irréductible apériodique et stationnaire, nous donnons une condition nécessaire et suffisante pour que les innovations $(V_n)_{n\in\mathbf{Z}}$ déterminent complètement la chaîne $(X_n)_{n\in\mathbf{Z}}$. Cette condition généralise et améliore la condition suffisante, mais non nécessaire, donnée par Rosenblatt dans [11].

Dans la troisième partie, nous nous plaçons dans le cadre d'une chaîne de Markov homogène sur un espace d'états fini et nous décrivons plus précisément l'information manquante lorsque la suite d'innovations ne détermine pas complètement la chaîne.

**1. Chaînes de Markov constructives inhomogènes.** Dans cette partie, $X = (X_n)_{n\in\mathbf{Z}}$ est une chaîne de Markov inhomogène gouvernée par une suite de variables $V = (V_n)_{n\in\mathbf{Z}}$ et une suite d'applications $(f_n)_{n\in\mathbf{Z}}$. Plus précisément, pour tout $n \in \mathbf{Z}$:

- $X_n$ est une variable aléatoire à valeurs dans un espace mesurable $(E_n, \mathcal{E}_n)$;
- $V_{n+1}$ est une variable aléatoire à valeurs dans un espace mesurable $(G_{n+1}, \mathcal{G}_{n+1})$, indépendante de la suite $((X_k, V_k))_{k\leq n}$;
- $X_{n+1} = f_n(X_n, V_{n+1})$ où $f_n$ est une application mesurable de $(E_n \times G_{n+1}, \mathcal{E}_n \otimes \mathcal{G}_{n+1})$ dans $(E_{n+1}, \mathcal{E}_{n+1})$.

Nous allons étudier la loi conditionnelle de $X$ sachant $\sigma(V) \vee \mathcal{F}_{-\infty}^{X,V}$. Mais voyons d'abord deux résultats simples sur la tribu asymptotique $\mathcal{F}_{-\infty}^{X,V}$.

1.1. *Tribu asymptotique d'une chaîne constructive inhomogène.* En supposant que toutes les tribus sont complétées, on a alors:

PROPOSITION 1.  *Les tribus asymptotiques $\mathcal{F}_{-\infty}^{(X,V)}$ et $\mathcal{F}_{-\infty}^{X}$ sont égales.*

PROOF. Pour montrer l'inclusion non triviale, il suffit de démontrer que pour tout $A \in \mathcal{F}_0^{(X,V)}$, $P[A|\mathcal{F}_{-\infty}^{(X,V)}] = P[A|\mathcal{F}_{-\infty}^{X}]$.

Par un argument de classes monotones, il suffit de se limiter au cas où $A \in \sigma(((X_k, V_k))_{N+1\leq k\leq 0})$ avec $N \in -\mathbf{N}^*$. Mais à cause de la relation de récurrence $X_{n+1} = f_n(X_n, V_{n+1})$, un tel événement $A$ s'exprime en fonction de $X_N, V_{N+1}, \ldots, V_0$. Plus généralement, pour tout $n \leq N$ il s'exprime en fonction de $X_n, V_{n+1}, \ldots, V_0$. Par indépendance de $(V_{n+1}, \ldots, V_0)$ et de $\mathcal{F}_n^{(X,V)}$, $P[A|\mathcal{F}_n^{(X,V)}]$ est une fonction borélienne de $X_n$. Donc $P[A|\mathcal{F}_{-\infty}^{(X,V)}] = \lim_{n\to-\infty} P[A|\mathcal{F}_n^{(X,V)}]$ est mesurable pour la tribu asymptotique $\mathcal{F}_{-\infty}^{X}$. □

Dans la suite, on notera simplement $\mathcal{F}_{-\infty}$ cette tribu. Comme cette tribu joue un rôle important par la suite, on peut se demander si elle est triviale (c'est-à-dire formée uniquement d'événements de probabilité 0 ou 1). La proposition qui suit donne une condition suffisante pour qu'elle le soit.

PROPOSITION 2. *Si $(X_n)_{n\in\mathbf{N}}$ est une chaîne de Markov homogène, stationnaire, irréductible et apériodique sur un espace d'états dénombrable, alors $\mathcal{F}_{-\infty}$ est triviale.*

PROOF. Notons $\Pi$ le noyau de transition de la chaîne et $\pi$ la probabilité invariante de $\Pi$. Il suffit de démontrer que pour tout $N \in -\mathbf{N}^*$ et pour tout $A \in \sigma(X_N, \ldots, X_0)$, $P[A|\mathcal{F}^X_{-\infty}] = P(A)$. Mais pour un tel événement $A$, $P[A|\mathcal{F}^X_N] = h(X_N)$ avec $h$ bornée de $E$ dans $\mathbf{R}$. Pour tout $n \leq N$, on a donc $P[A|\mathcal{F}^X_n] = \Pi_{N-n}h(X_n)$. Comme $X_n$ a pour loi $\pi$, on a donc

$$\|P[A|\mathcal{F}^X_n] - P(A)\|_1 = \int_E |\Pi^{N-n}h(x) - P(A)|\pi(dx),$$

mais pour tout $x \in E$, $\Pi^{N-n}(x, \cdot)$ converge en variation totale vers $\pi$ quand $n \to -\infty$, d'où $\Pi^{N-n}h(x) \to \int_E h\,d\pi = \mathbf{E}[h(X_N)] = P(A)$. Par convergence dominée, on a donc $\|P[A|\mathcal{F}^X_n] - P(A)\|_1 \to 0$ quand $n \to -\infty$, ce qui montre que $P[A|\mathcal{F}^X_{-\infty}] = P(A)$ par le théorème de convergence des martingales. □

REMARQUE. La démonstration et le résultat ci-dessus restent valables pour un espace d'état quelconque à condition de supposer que la chaîne est récurrente positive au sens de Harris (voir [10]).

Signalons également que la stationnarité est essentielle. Pour s'en convaincre, il suffit de considérer les probabilités de transition sur $\mathbf{N}$ donnés par $p(x, x-1) = 1$ pour tout $x \in \mathbf{N}^*$ et $p(0, y) = (\frac{1}{2})^{y+1}$ pour tout $y \in \mathbf{N}$. La matrice de transition est irréductible, apériodique et la probabilité donnée par $\pi\{y\} = (\frac{1}{2})^{y+1}$ est invariante. On peut construire une chaîne de Markov indexée par $\mathbf{Z}$ en prenant une variable aléatoire $L$ à valeurs dans $\mathbf{N}$, et en posant $X_n = L - n$ pour tout $n \in -\mathbf{N}$. La variable aléatoire $L$ est alors mesurable pour la tribu asymptotique $\mathcal{F}_{-\infty}$.

1.2. *La loi conditionnelle de $X$ connaissant $V$.* Nous supposons que les espaces d'états $(E_n, \mathcal{E}_n)$ sont lusiniens (espaces de Borel suivant la terminologie de Kallenberg [5]). Il en est alors de même pour l'espace $\prod_{k\in\mathbf{Z}} E_k$ muni de la tribu produit, ce qui permet de parler de la loi conditionnelle de $X$ sachant $\sigma(V) \vee \mathcal{F}_{-\infty}$. Nous avons le résultat suivant:

THÉORÈME 3. *Presque sûrement, la loi de $X$ sachant $\sigma(V) \vee \mathcal{F}_{-\infty}$ est diffuse ou uniforme sur un ensemble fini. De plus, le nombre d'atomes (et donc leur masse) est une variable aléatoire $\mathcal{F}_{-\infty}$-mesurable.*

De ce théorème, on déduit immédiatement:

COROLLAIRE 4. *Si la tribu $\mathcal{F}_{-\infty}$ est triviale, il existe une variable aléatoire $U$, indépendante de $V$ et de loi uniforme sur $[0,1]$ ou sur un ensemble fini, telle que $\sigma(X) = \sigma(V) \vee \sigma(U)$.*

La démonstration de ces résultats utilise fortement des résultats classiques sur les lois conditionnelles que nous rappelons ci-dessous.

LEMME 5. *Soient $X$ une variable aléatoire définie sur un espace probabilisé $(\Omega, \mathcal{A}, P)$ et à valeurs dans un espace de Borel $(E, \mathcal{E})$, $\mathcal{F}$ et $\mathcal{G}$ deux sous-tribus de $\mathcal{A}$ et $(\nu_\omega)_{\omega \in \Omega}$ une version régulière de la loi de $X$ sachant $\mathcal{F}$. Alors:*

*1. Quelle que soit l'application mesurable $f$ de $(\Omega \times E, \mathcal{G} \otimes \mathcal{E})$ dans $\mathbf{R}_+$, l'application $\omega \mapsto \int_E f(\omega, x)\nu_\omega(dx)$ est $\mathcal{F} \vee \mathcal{G}$-mesurable.*

*2. Si $f$ est mesurable de $(\Omega \times E, \mathcal{F} \otimes \mathcal{E})$ dans $\mathbf{R}_+$, la variable aléatoire $\omega \mapsto \int_E f(\omega, x)\nu_\omega(dx)$ est l'espérance conditionnelle de $\omega \mapsto f(\omega, X(\omega))$ sachant $\mathcal{F}$.*

*3. Plus généralement, si $f$ est une application mesurable de $(\Omega \times E, \mathcal{F} \otimes \mathcal{E})$ dans un espace mesurable $(E', \mathcal{E}')$. Alors les mesures images des mesures $\nu_\omega$ par les applications $x \mapsto f(\omega, x)$ fournissent une version régulière de la loi conditionnelle de $\omega \mapsto f(\omega, X(\omega))$ sachant $\mathcal{F}$.*

Le point 1 est une application directe du théorème des classes monotones. Le point 2 fait l'objet du Théorème 5.4 de Kallenberg [5]. Le point 3 s'en déduit immédiatement.

Passons maintenant à la démonstration du théorème.

PROOF OF THEOREM 3. Soit $\mathcal{F} = \sigma(V) \vee \mathcal{F}_{-\infty}$ et soit $(\nu_\omega)_{\omega \in \Omega}$ une version régulière de la loi de $X$ sachant $\mathcal{F}$. Notons $A(\omega) = \nu_\omega\{X(\omega)\} = \int \mathbf{I}_{\{X(\omega)=x\}}\nu_\omega(dx)$. Le point 1 du Lemme 5 montre que $A$ est une variable aléatoire mesurable pour $\sigma(X, V)$. Nous allons voir que $A$ est même $\mathcal{F}_{-\infty}$-mesurable.

Pour $n \in \mathbf{Z}$ et $\omega \in \Omega$, notons $\nu_n(\omega, \cdot)$ la mesure image de $\nu_\omega$ par la projection canonique $\pi_n : x = (x_k)_{k \in \mathbf{Z}} \mapsto x_{n]} = (x_k)_{k \leq n}$ de $\prod_{k \in \mathbf{Z}} E_k$ dans $\prod_{k \leq n} E_k$. La famille $(\nu_n(\omega, \cdot))_{\omega \in \Omega}$ est une version régulière de la loi de $X_{n]}$ sachant $\mathcal{F}$.

Mais comme $X$ est une fonction de $X_{n]}$ et de $V$, disons $X = F_n(X_{n]}, V)$, le point 3 du Lemme 5 nous assure que pour presque tout $\omega \in \Omega$, la loi $\nu_\omega$ est l'image de la loi $\nu_n(\omega, \cdot)$ par l'application $F_n(\cdot, V(\omega))$. Par conséquent, les atomes de $\nu_n(\omega, \cdot)$ sont les images par $\pi_n$ des atomes de $\nu_\omega$ et il y a conservation des masses. Donc pour presque tout $\omega \in \Omega$, $A(\omega) = \nu_\omega\{X(\omega)\} = \nu_n(\omega, \{X_{n]}(\omega)\})$.

Par indépendance de $\mathcal{F}_n^{(X,V)} = \sigma(X_{n]}) \vee \sigma(V_{n]}) \vee \mathcal{F}_{-\infty}$ et de $V_{[n+1} = (V_k)_{k \geq n+1}$, $\nu_n$ est la loi conditionnelle de $X_{n]}$ sachant $\sigma(V_{n]}) \vee \mathcal{F}_{-\infty}$. En particulier $\nu_n$ est une mesure aléatoire $\mathcal{F}_n^{(X,V)}$-mesurable. Donc $A = \nu_n(\{X_{n]}\})$ est une variable aléatoire $\mathcal{F}_n^{(X,V)}$-mesurable pour tout $n \in \mathbf{Z}$, ce qui montre que $A$ est $\mathcal{F}_{-\infty}$-mesurable.

Dans l'égalité $\mathbf{I}_{[\nu_\omega\{X(\omega)\}=A(\omega)]} = 1$, prenons les espérances conditionnelles sachant $\mathcal{F}$ en appliquant le point 2 du Lemme 5 à l'espace $(E,\mathcal{E}) = (\prod_{k \in \mathbf{Z}} E_k, \bigotimes_{k \in \mathbf{Z}} \mathcal{E}_k)$ et à l'application $f : (\omega, x) \mapsto \mathbf{I}_{[\nu_\omega\{x\}=A(\omega)]}$. On obtient pour presque tout $\omega \in \Omega$,

$$\int_E \mathbf{I}_{[\nu_\omega\{x\}=A(\omega)]} \nu_\omega(dx) = 1.$$

Mais l'intégrale du membre de gauche vaut $A(\omega)$ fois le nombre d'atomes de masse $A(\omega)$ de $\nu_\omega$ si $A(\omega) > 0$ et est égale à la masse de la partie diffuse de $\nu_\omega$ si $A(\omega) = 0$. La mesure $\nu_\omega$ est donc uniforme sur un ensemble fini de cardinal $1/A(\omega)$ si $A(\omega) > 0$ et diffuse si $A(\omega) = 0$. □

REMARQUE. La démonstration montre également que presque sûrement, la loi de $X_{n]}$ sachant $\sigma(V_{n]}) \vee \mathcal{F}_{-\infty}$ est diffuse ou uniforme sur un ensemble fini, avec un nombre d'atomes qui ne dépend pas de $n$.

Pour démontrer le Corollaire 4 à partir du Théorème 3 on utilise à nouveau le point 3 du Lemme 5, en choisissant une application $f$ mesurable de $(\Omega \times E, \mathcal{F} \otimes \mathcal{E})$ dans $\mathbf{R}$ telle que la mesure image de $\nu_\omega$ par l'application $f(\cdot, \omega)$ soit la loi uniforme sur $[0,1]$ si $\nu_\omega$ est diffuse et loi uniforme sur $[1, \ldots, N]$ si $\nu_\omega$ est uniforme sur un ensemble de cardinal $N$. Si la tribu $\mathcal{F}_{-\infty}$ est triviale, le nombre $N$ d'atomes de $\nu_\omega$ est presque sûrement constant. La variable aléatoire $U : \omega \mapsto f(X(\omega), \omega)$ est alors indépendante de $\mathcal{F}$ puisque sa loi conditionnelle sachant $\mathcal{F}$ est presque sûrement constante (égale à la loi uniforme sur $[0,1]$ si $N = 0$, sur $[1, \ldots, N]$ si $N \geq 1$).

1.3. *Exemple*: *marche aléatoire sous l'action d'un groupe compact.* Considérons un groupe compact $G$ agissant continûment sur un espace polonais $(E, d)$ et une probabilité $\mu$ dont le support n'est contenu dans aucune classe modulo un sous-groupe distingué non trivial.

Soit $(X_n)_{n \in \mathbf{Z}}$ une marche aléatoire à valeurs dans $E$ gouvernée par une suite $(V_n)_{n \in \mathbf{Z}}$ de variables aléatoires indépendantes de loi $\mu$ sur $G$, c'est-à-dire telle que pour tout $n \in \mathbf{Z}$, $X_{n+1} = V_{n+1} \cdot X_n$ avec $V_{n+1}$ indépendante de la suite $((X_k, V_k))_{k \leq n}$.

Notons $\pi$ la projection canonique de $E$ sur l'espace $G \setminus E$ des orbites sous l'action de $G$. Munissons $G \setminus E$ de la tribu quotient, c'est-à-dire de la plus grande tribu rendant $\pi$ mesurable.

THÉORÈME 6. *Sous les conditions, précédentes:*

- La tribu $\mathcal{F}_{-\infty}$ est engendrée par la variable aléatoire $\pi(X_0) = G \cdot X_0$ (à valeurs dans $G \setminus E$).
- La variable aléatoire $X_0$ est indépendante de $V$ conditionnellement à $\mathcal{F}_{-\infty}$.
- La de loi $X_0$ sachant $\mathcal{F}_{-\infty}$ est la loi uniforme sur $G \cdot X_0$, c'est-à-dire l'image de la mesure de Haar sur $G$ par l'application $g \mapsto g \cdot X_0$.

REMARQUE PRÉLIMINAIRE. La variable aléatoire $\pi(X_0)$ est $\mathcal{F}_{-\infty}^X$-mesurable car $\pi(X_n) = \pi(X_0)$ pour tout $n \in \mathbf{Z}$.

Introduisons une notation: si $\phi$ est une application mesurable bornée de $E$ dans $\mathbf{R}$, on définit une application $\bar{\phi}$ mesurable bornée de $G \setminus E$ dans $\mathbf{R}$ en prenant la moyenne orbite par orbite. Autrement dit, pour tout $x \in E$,

$$\bar{\phi}(G \cdot x) = \int_G \phi(g \cdot x)\, dg,$$

où $dg$ est la mesure de Haar sur $G$. Avec cette notation, nous pouvons énoncer le lemme sur lequel repose le théorème:

LEMME 7. *Soient $\phi$ continue bornée de $E$ dans $\mathbf{R}$ et $\psi$ continue bornée de $G^{N+1}$ dans $\mathbf{R}$. Alors*

$$\mathbf{E}[\phi(X_0)\psi(V_0,\ldots,V_{-N})|\mathcal{F}_{-\infty}^{(X,V)}] = \bar{\phi}(\pi(X_0))\mathbf{E}[\psi(V_0,\ldots,V_{-N})].$$

PROOF. Quel que soit $n \in \mathbf{N}^*$,

$$X_0 = V_0 \cdots V_{-N} V_{-N-1} \cdots V_{-N-n+1} \cdot X_{-N-n}.$$

Par indépendance de $(V_0,\ldots,V_{-N-n+1})$ et de $\mathcal{F}_{-N-n}^{(X,V)}$ on a donc

$$\mathbf{E}[\phi(X_0)\psi(V_0,\ldots,V_{-N})|\mathcal{F}_{-N-n}^X] = h_n(X_{-N-n}),$$

où $h_n(x) = \mathbf{E}[\phi(V_0 \cdots V_{-N} V_{-N-1} \cdots V_{-N-n+1} \cdot x)\psi(V_0,\ldots,V_{-N})]$. D'après les hypothèses faites sur $\mu$, la loi de $V_{-N-1} \cdots V_{-N-n+1}$ converge étroitement vers la mesure de Haar de $G$ (voir [8], Théorème 4.8, ou [4], Théorème 2.1.4), donc pour tout $(v_0,\ldots,v_{-N}) \in G^{N+1}$,

$$\mathbf{E}[\phi(v_0 \cdots v_{-N} V_{-N-1} \cdots V_{-N-n+1} \cdot x)] \to \int_G \phi(v_0 \cdots v_{-N} g \cdot x) = \bar{\phi}(\pi(x)).$$

En conditionnant l'expression de $h_n(x)$ par rapport à $(V_0,\ldots,V_{-N})$, on en déduit que pour tout $x \in E$,

$$h_n(x) \to \bar{\phi}(\pi(x))\mathbf{E}[\psi(V_0,\ldots,V_{-N})].$$

Mais par compacité de $G$, la restriction de $\phi$ à chaque orbite est uniformément continue. Par conséquent, la famille des restrictions des $h_n$ à chaque orbite est équicontinue. Comme la suite $(X_{-N-n})_{n \in \mathbf{N}^*}$ reste sur l'orbite $\pi(X_0)$, on a donc

$$h_n(X_{-N-n}) \to \bar{\phi}(\pi(X_0))\mathbf{E}[\psi(V_0, \ldots, V_{-N})],$$

ce qui montre le lemme par application du théorème de convergence des martingales. $\square$

Passons maintenant à la démonstration du Théorème 6.

PROOF OF THEOREM 6. Par un argument de classes monotone, le Lemme 5 montre que pour toute variable aléatoire bornée $Z$ mesurable pour la tribu $\mathcal{F}_0^{(X,V)} = \sigma(X_0) \vee \mathcal{F}_0^V$, la variable aléatoire $\mathbf{E}[Z|\mathcal{F}_{-\infty}^{(X,V)}]$ est une fonction mesurable de $\pi(X_0)$. D'où le premier point, compte tenu de la remarque préliminaire.

Le lemme montre également que conditionnellement à $\mathcal{F}_{-\infty}$, $X_0$ et $\mathcal{F}_0^V$ sont indépendantes, donc $X_0$ et $\mathcal{F}_\infty^V$ sont indépendantes, puisque la suite $(V_n)_{n \in \mathbf{N}^*}$ est indépendante de $\mathcal{F}_0^{(X,V)}$.

Le dernier point du théorème découle du lemme en prenant $\psi$ constante égale à 1. $\square$

**2. Une chaîne stationnaire sur un espace dénombrable est-elle déterminée par ses innovations?** Dans cette partie, on se donne un ensemble dénombrable $E$, un espace mesurable $(G, \mathcal{G})$, une application mesurable $f$ de $E \times G$ dans $E$ et une probabilité $\beta$ sur $(G, \mathcal{G})$. La question qui nous intéresse ici est de savoir à quelle condition une chaîne de Markov stationnaire gouvernée par la relation de récurrence $X_{n+1} = f(X_n, V_{n+1})$, avec $(V_n)_{n \in \mathbf{Z}}$ indépendantes et de loi $\beta$, est déterminée par la suite d'innovations $(V_n)_{n \in \mathbf{Z}}$.

2.1. *La condition suffisante de Rosenblatt.* En introduisant les applications $f_v : x \mapsto f(x, v)$ de $E$ dans $E$, on peut réécrire la relation de récurrence sous la forme $X_{n+1} = f_{V_{n+1}}(X_n)$, ce qui donne par récurrence $X_n = f_{V_n} \circ \cdots \circ f_{V_{m+1}}(X_m)$ pour tout $m \leq n$. Cette réécriture fournit une condition suffisante pour que la suite $(V_n)_{n \in \mathbf{Z}}$ détermine complètement la chaîne $(X_n)_{n \in \mathbf{Z}}$: il suffit qu'il existe $l \in \mathbf{N}^*$ et $c \in E$ tels que la probabilité que la composée $f_{V_l} \circ \cdots \circ f_{V_1}$ soit constante égale à $c$ avec probabilité strictement positive.

Dans ce cas en effet, presque sûrement, pour tout $n \in \mathbf{Z}$, on peut trouver $m \leq n$ tel que $f_{V_m} \circ \cdots \circ f_{V_{m-l+1}}$ soit constante égale à $c$. Mais pour tout $m \leq n$, l'événement $[f_{V_m} \circ \cdots \circ f_{V_{m-l+1}} = c]$ appartient à $\mathcal{F}_n^V$, et sur cet événement, $X_m = c$, d'où $X_n = f_{V_n} \circ \cdots \circ f_{V_{m+1}}(c)$. On en déduit que $X_n$ est (presque sûrement) une fonction de $(V_k)_{k \leq n}$, autrement dit que $X_n$ est mesurable pour $\mathcal{F}_n^V$. On peut remarquer que la stationnarité de la chaîne $(X_n)_{n \in \mathbf{Z}}$ n'a pas servi dans la démonstration, mais on vérifie facilement que la stationnarité découle de l'hypotyèse que $f_{V_l} \circ \cdots \circ f_{V_1}$ soit constante égale à $c$ avec probabilité strictement positive.

Cette condition est l'ingrédient essentiel de l'algorithme coupling from the past de Propp et Wilson [9], permettant une simulation exacte le la loi stationnaire d'une chaîne de Markov récurrente positive. Mais on trouve déjà une condition semblable obtenue de cette façon (dans un cas particulier) dans un article de Rosenblatt de 1959: dans le Paragraphe 3 de [11], Rosenblatt s'intéresse aux chaînes de Markov uniformes sur un espace d'états dénombrable, c'est-à-dire aux chaînes de Markov dont les probabilités de transition issues d'un état et rangées par ordre décroissant ne dépendent pas de l'état considéré (autrement dit, les lignes de la matrice de transition se déduisent les unes des autres par permutation de leur coefficients). Cette hypothèse permet de construire des innovations à partir de la chaîne $(X_n)_{n\in\mathbf{N}}$, sans apport d'aléa extérieur.

Plus précisément, la loi de $X_{n+1}$ sachant que $X_n = x$ est atomique et les masses rangées par ordre décroissant $p_1 \geq p_2 \geq \cdots$ ne dépendent pas de $x$. Des bijections $\phi_x$ convenables de $E$ dans $\{1, 2, \ldots\}$ transforment ces lois en la loi $p_1\delta_1 + p_2\delta_2 + \cdots$, si bien que la variable aléatoire $\xi_{n+1} = \phi_{X_n}(X_{n+1})$ est indépendante de $\mathcal{F}_n^X$ et de loi $p_1\delta_1 + p_2\delta_2 + \cdots$. De plus, on peut inverser la formule définissant $\xi_{n+1}$ et la réécrire sous la forme $X_{n+1} = \phi_{X_n}^{-1}(\xi_{n+1}) = f(X_n, \xi_{n+1})$.

Dans son Théorème 1, Rosenblatt affirme qu'une condition nécessaire et suffisante pour que les innovations $(\xi_n)_{n\in\mathbf{N}}$ déterminent la chaîne $(X_n)_{n\in\mathbf{N}}$ est que le semigroupe engendré par les applications $f_k : x \mapsto f(x, k)$ soit point collapsing, c'est-à-dire qu'une composée convenable des applications $f_k$ soit constante. En réalité, Rosenblatt regarde le semigroupe engendré par les matrices $M_k = (\mathbf{I}_{[j=f(i,k)]})_{i,j\in E}$, mais cela revient au même.

Cependant, Rosenblatt ne démontre vraiment qu'une implication: si une composée d'applications $f_k$ est constante, alors la suite $(\xi_n)_{n\in\mathbf{Z}}$ engendre la même filtration que la chaîne de Markov $(X_n)_{n\in\mathbf{Z}}$. Mais pour la réciproque, l'argument de Rosenblatt revient à admettre l'idée intuitive suivante: si l'on connaît la suite $(\xi_n)_{n\in\mathbf{Z}}$, la seule chose qu'on puisse dire de $X_n$ est que $X_n$ appartient à l'intersection des images des composées $f_{\xi_n} \circ f_{\xi_{n-1}} \circ \cdots \circ f_{\xi_{n-l+1}}$ pour $l \in \mathbf{N}$. Cet argument heuristique est juste dans le cas où $E$ est fini, et nous en donnerons un énoncé rigoureux dans le Théorème 11. Mais il est faux lorsque $E$ est infini. Signalons que l'équivalence dans le cas où $E$ est fini est démontrée par Laurent dans sa thèse de doctorat [6] (Proposition 9.3.3). La preuve pour le sens difficile utilise un argument de couplage (couplage de Doeblin) différent du nôtre.

Nous allons établir dans le paragraphe suivant une condition nécessaire et suffisante. Ce résultat montre que la condition de Rosenblatt est effectivement nécessaire lorsque $E$ est fini, mais qu'elle ne l'est plus dans le cas général, même pour une chaîne uniforme.

**2.2. Une condition nécessaire et suffisante.** Soit $(V_n)_{n\in\mathbf{Z}}$ une suite de variables aléatoires indépendantes de loi $\beta$.

Pour tout $x \in E$, notons $\Pi(x,\cdot)$ la loi de $f(x,V_1)$. Alors $\Pi$ est le noyau de toute chaîne de Markov gouvernée par la relation de récurrence $X_{n+1} = f(X_n, V_{n+1})$.

Pour tout $(x,y) \in E^2$, notons $\Pi_2((x,y),\cdot)$ la loi de $(f(x,V_1), f(y,V_1))$. De même $\Pi_2$ est le noyau de toute chaîne de Markov sur $E^2$ gouvernée par la relation de récurrence $(X_{n+1}, Y_{n+1}) = (f(X_n, V_{n+1}), f(Y_n, V_{n+1}))$.

Notons $D$ la diagonale de $E^2$ et définissons une relation sur $E$, que nous appelerons accordabilité.

DÉFINITION 8. On dit que deux états $x$ et $y$ de $E$ sont accordables si $D$ est accessible depuis $(x,y)$ pour les chaînes de noyau $\Pi_2$ (autrement dit, s'il existe $l \in \mathbf{N}^*$ tel que $P[f_{V_l} \circ \cdots \circ f_{V_1}(x) = f_{V_l} \cdots \circ f_{V_1}(y)] > 0$).

REMARQUE. Cette définition n'est pas affectée si l'on remplace la loi $\beta$ par une probabilité équivalente (ou même par une probabilité pour laquelle la loi de $f_{V_1}$ est changée en une probabilité équivalente). Les résultats qui suivent ne dépendent donc que de la classe d'équivalence de la loi $\beta$.

Énonçons un résultat pour les chaînes récurrentes positives.

THÉORÈME 9. *Si le noyau $\Pi$ est irréductible et positivement récurrent, les conditions suivantes sont équivalentes:*

(1) *Les points de $E$ sont deux-à-deux accordables.*
(2) *Pour toute chaîne stationnaire $(X_n)_{n\in\mathbf{Z}}$ gouvernée par $f$ et par une suite $(V_n)_{n\in\mathbf{Z}}$ de variables aléatoires indépendantes de loi $\beta$, on a $\mathcal{F}^X_\cdot \subset \mathcal{F}^V_\cdot$.*
(3) *Toute chaîne stationnaire pour le noyau $\Pi_2$ vit sur la diagonale $D$.*

PROOF. Nous montrons les équivalences $(1) \Leftrightarrow (3)$ et $(2) \Leftrightarrow (3)$.

$(1) \Rightarrow (3)$ Soit $(X_n, Y_n)_{n\in\mathbf{Z}}$ une chaîne stationnaire pour le noyau $\Pi_2$. Soient $x$ et $y$ deux points distincts de $E$. Comme $D$ est un ensemble absorbant pour $\Pi_2$ et accessible depuis $(x,y)$, l'état $(x,y)$ est transient. Donc $P[(X_n, Y_n) = (x,y)] \to 0$ quand $n \to \infty$, d'où $P[(X_n, Y_n) = (x,y)] = 0$ par stationnarité. Cela montre que la chaîne $(X_n, Y_n)_{n\in\mathbf{Z}}$ vit sur $D$.

Non $(1) \Rightarrow$ non $(3)$. Soit $\pi$ la probabilité invariante pour $\Pi$. Supposons que $a$ et $b$ sont deux états non accordables de $E$. Nous allons construire une probabilité invariante pour $\Pi_2$ portée par $E^2 \setminus D$, ce qui contredira $(3)$. Pour tout $n \in \mathbf{N}$, notons $\nu_n = \delta_{(a,b)}\Pi_2^n$ la loi de la position à l'instant $n$ d'une chaîne issue de $(a,b)$. Les lois marginales de $\nu_n$ sont $\delta_a\Pi^n$ et $\delta_b\Pi^n$, qui sont majorées par les mesures finies $\pi/\pi\{a\}$ et $\pi/\pi\{b\}$. Pour toute partie finie $F \subset E$, on a donc

$$\nu_n(E^2 \setminus F^2) \leq \pi(E \setminus F)\Big(\frac{1}{\pi\{a\}} + \frac{1}{\pi\{b\}}\Big).$$

Il s'en suit que la suite $(\nu_n)_{n\geq 0}$ est tendue, ainsi que la suite de ses moyennes de Cesàro $((\nu_0+\cdots+\nu_{n-1})/n)_{n\geq 1}$. Un argument classique montre que toute valeur d'adhérence de cette suite est invariante pour $\Pi_2$. Mais $\nu_n(D) = 0$ pour tout $n \in \mathbf{N}$ puisque $a$ et $b$ sont non accordables, d'où le résultat.

$(2) \Rightarrow (3)$ Soit $(X_n, Y_n)_{n \in \mathbf{Z}}$ est une chaîne stationnaire pour le noyau $\Pi_2$. Notons $\alpha$ la loi de $(X_0, Y_0)$.

Montrons d'abord qu'on peut construire, sur un espace probabilisé convenable, une chaîne $(X'_n, Y'_n)_{n \in \mathbf{Z}}$ de même loi que $(X_n, Y_n)_{n \in \mathbf{Z}}$ et gouvernée par une suite d'innovations $(V'_n)_{n \in \mathbf{Z}}$ de même loi que $(V_n)_{n \in \mathbf{Z}}$. Remarquons que si une telle chaîne existe, le processus $(X'_n, Y'_n, V'_{n+1})_{n \in \mathbf{Z}}$ est une chaîne de Markov pour le noyau $\widetilde{\Pi}$ où $\widetilde{\Pi}((x,y,v), \cdot) = \delta_{(f(x,v), f(y,v))} \otimes \beta$. L'existence d'un telle chaîne découle du fait que la loi $\alpha \otimes \beta$ est invariante pour le noyau $\widetilde{\Pi}$.

Pour alléger les notations, on notera simplement $(X_n)_{n \in \mathbf{Z}}$, $(Y_n)_{n \in \mathbf{Z}}$ et $(V_n)_{n \in \mathbf{Z}}$ les processus ainsi construits. Les processus $X$ et $Y$ sont des chaînes de Markov stationnaires gouvernées par les relations de récurrence $X_{n+1} = f(X_n, V_{n+1})$ $Y_{n+1} = f(Y_n, V_{n+1})$. D'après l'hypothèse, pour tout $n \in \mathbf{Z}$, il existe deux applications $F_1$ et $F_2$ mesurables de $G^{]-\infty,\ldots,n]}$ dans $E$, telle que $X_n = F_1(V_{n]})$ et $Y_n = F_2(V_{n]})$ presque sûrement.

Pour montrer que $X_n = Y_n$ presque sûrement, il suffit alors de vérifier que $(X_n, V_{n]})$ et $(Y_n, V_{n]})$ ont même loi. Pour cela, on remarque que pour tout $m < n$, $(X_m, V_{m+1}, \ldots, V_n)$ et $(Y_m, V_{m+1}, \ldots, V_n)$ ont même loi $\pi \otimes \beta^{n-m}$. Comme $X_n = f_{V_n} \circ \cdots \circ f_{V_{m+1}}(X_m)$ et $Y_n = f_{V_n} \circ \cdots \circ f_{V_{m+1}}(Y_m)$, les variables aléatoires $(X_n, V_{m+1}, \ldots, V_n)$ et $(Y_n, V_{m+1}, \ldots, V_n)$ ont aussi même loi, ce qui montre le résultat.

$(3) \Rightarrow (2)$ Soient $V = (V_n)_{n \in \mathbf{Z}}$ une suite de variables aléatoires indépendantes de loi $\beta$ et $X = (X_n)_{n \in \mathbf{Z}}$ une chaîne de Markov gouvernée par $f$ et par la suite $V$. Nous allons montrer que la loi conditionnelle de $X$ connaissant $V$ est presque sûrement une masse de Dirac. Pour cela, considérons (en agrandissant si besoin l'espace probabilisé) une variable aléatoire $Y = (Y_n)_{n \in \mathbf{Z}}$ à valeurs dans $E^{\mathbf{Z}}$ telle que conditionnellement à $V$:

- $Y$ est indépendante de $X$;
- $Y$ a même loi que $X$.

La suite $(X_n, Y_n)_{n \in \mathbf{Z}}$ est alors une chaîne stationnaire pour le noyau $\Pi_2$. Donc, presque sûrement, $X = Y$ d'où $P[X = Y | \sigma(V)] = 1$, ce qui montre que la loi conditionnelle de $X$ connaissant $V$ est une masse de Dirac. On a donc $\mathcal{F}^X_{+\infty} \subset \mathcal{F}^V_{+\infty}$. Pour montrer que $\mathcal{F}^X_n \subset \mathcal{F}^V_n$, il suffit de remarquer que pour tout $A \in \mathcal{F}^X_n$, on a par indépendance de $(V_k)_{k \geq n+1}$ et de $\mathcal{F}^X_n \vee \mathcal{F}^V_n$,

$$P[A|\mathcal{F}^V_n] = P[A|\mathcal{F}^V_{+\infty}] = P[A|\mathcal{F}^X_{+\infty}] = \mathbf{I}_A,$$

car $A \in \mathcal{F}^X_{+\infty} \subset \mathcal{F}^V_{+\infty}$, d'où $A \in \mathcal{F}^V_n$. $\square$

Remarques.

- Dans l'énoncé (2), on peut remplacer l'inclusion des filtrations $\mathcal{F}_\cdot^X \subset \mathcal{F}_\cdot^V$ par la condition plus faible $\mathcal{F}_{+\infty}^X \subset \mathcal{F}_{+\infty}^V$.
- L'hypothèse de récurrence positive ne sert qu'à assurer l'existence de $\pi$ dans la démonstration de l'implication $(3) \Rightarrow (1)$.
- Les énoncés équivalents (1), (2) et (3) entraînent que le noyau $\Pi$ est apériodique.

Le Théorème 9 montre que la condition de Rosenblatt est nécessaire lorsque $E$ est fini. En effet, une récurrence immédiate montre que l'accordabilité deux-à-deux entraîne l'existence d'un entier $l$ tel que la composée $f_{V_l} \circ \cdots \circ f_{V_1}$ soit constante avec probabilité $> 0$. Mais le Théorème 9 permet aussi de voir que la condition de Rosenblatt n'est pas nécessaire lorsque $E$ est infini, même pour une chaîne de Markov uniforme, comme le montre le contre-exemple ci-dessous.

Contre-exemple. Posons $E = \mathbf{N}$, $G = \{1, 2\}$, $f_1(x) = f(x, 1) = (x - 1)_+$ et $f_2(x) = f(x, 2) = x + 1$ et prenons une suite $(V_n)_{n \in \mathbf{Z}}$ de variables aléatoires indépendantes de loi $\frac{2}{3}\delta_1 + \frac{1}{3}\delta_2$. Notons $\Pi(x, \cdot)$ la loi de $f(x, V_1)$. Le noyau $\Pi$ ainsi défini est irréductible, apériodique et la probabilité $\pi$ donnée par $\pi\{y\} = (\frac{1}{2})^{y+1}$ pour tout $y \in \mathbf{N}$ est invariante.

Comme $f_1 \circ f_2 = \mathrm{id}_E$, une composée d'applications $f_1$ et $f_2$ s'écrit toujours $f_2^q \circ f_1^p$, et on voit immédiatement qu'une telle composée n'est jamais constante. Cependant, pour tout $l \in \mathbf{N}^*$, la composée $f_1^l$ vaut 0 sur $[0, \ldots, l]$, ce qui montre que les états sont deux-à-deux accordables.

2.3. *Cas des chaînes non stationnaires.* Supposons que les états sont sont deux-à-deux accordables. Si $(X_n)_{n \in \mathbf{Z}}$ est une chaîne de Markov gouvernée par $f$ et par la suite d'innovations $(V_n)_{n \in \mathbf{Z}}$, peut-on encore dire que $\mathcal{F}_\cdot^X \subset \mathcal{F}_\cdot^V$ si $(X_n)_{n \in \mathbf{Z}}$ n'est pas stationnaire? Nous allons voir que non, même si le noyau $\Pi$ est irréductible et positivement récurrent.

Une première obstruction évidente se produit lorsque la tribu asymptotique $\mathcal{F}_{-\infty} = \mathcal{F}_{-\infty}^X = \mathcal{F}_{-\infty}^{X,V}$ n'est pas triviale, comme dans l'exemple du Paragraphe 1.1. Cependant, ce n'est pas la seule obstruction, et l'exemple ci-dessous montre qu'on peut avoir $\mathcal{F}_\cdot^X \not\subset \mathcal{F}_\cdot^V$ même si $\mathcal{F}_{-\infty}$ est triviale.

Exemple. Prenons $E = G = \mathbf{Z}$, soit $\gamma$ une probabilité sur $\mathbf{N}^*$ qui charge tous les points. Sur un espace probabilisé convenable, définissons deux suites de variables aléatoires toutes indépendantes: $(W_n)_{n \in \mathbf{Z}}$ de loi $\gamma$ et $(\varepsilon_n)_{n \in \mathbf{Z}}$ de loi $\frac{1}{2}(\delta_{-1} + \delta_1)$. Alors les variables aléatoires $V_n = \varepsilon_{n-1}\varepsilon_n W_n$ sont indépendantes et de même loi. De plus, leur loi commune $\beta$ charge tous les points de $\mathbf{Z}^*$.

Soit $f : E \times G \to E$ définie par $f(0, v) = v$ et $f(x, v) = (|x| - 1)\,\mathrm{sgn}(xv)$ si $x \neq 0$. On vérifie facilement que la suite $(X_n)_{n \in \mathbf{Z}}$ définie par $X_n = |n|\varepsilon_n$ si

$n \leq 0$ et $X_n = f_{V_n} \circ \cdots \circ f_{V_1}(0)$ si $n > 0$ est une chaîne de Markov gouvernée par $f$ et par la suite d'innovations $(V_n)_{n \in \mathbf{Z}}$.

La tribu $\mathcal{F}_{-\infty}$ est triviale puisque c'est la tribu asymptotique de la suite $(\varepsilon_n)_{n \in \mathbf{Z}}$. Le noyau $\Pi$ est irréductible, apériodique, récurrent et même positivement récurrent si $\mathbf{E}[W_1] < +\infty$. Les états sont accordables deux-à-deux. Pourtant, les innovations ne déterminent pas la chaîne puisque $X_{-1} = \varepsilon_{-1}$ est indépendante de $(V_n)_{n \in \mathbf{Z}}$.

**3. Description de l'information manquante pour une chaîne homogène sur un espace d'états fini.** Dans cette partie, $(X_n)_{n \in \mathbf{Z}}$ est une chaîne de Markov sur un espace d'états fini $E$, gouvernée par une suite $(V_n)_{n \in \mathbf{Z}}$ d'innovations i.i.d. à valeurs un espace d'états mesurable $(G, \mathcal{G})$ et par l'application mesurable $f$ de $E \times G$ dans $G$. On suppose que le noyau $\Pi$ est irréductible et apériodique, ce qui entraîne que la chaîne $(X_n)_{n \in \mathbf{Z}}$ est stationnaire et que la tribu $\mathcal{F}_{-\infty}$ est triviale. On note $\beta$ la loi commune des innovations $V_n$.

3.1. *Description de l'information manquante.* Nous allons décrire l'information fournie par la suite d'innovations $(V_n)_{n \in \mathbf{Z}}$, et montrer que l'information manquante est liée au nombre maximal $M$ d'états deux-à-deux non accordables. Le principal outil de la démonstration est l'étude d'une marche aléatoire associée à la suite d'applications $f_{V_n}$.

Soit $H$ l'image essentielle de $f_{V_1}$. Comme l'ensemble des applications de $E$ dans $E$ est fini, on a simplement $H = \{h \in E^E : P[f_{V_1} = h] > 0\}$. Soit $S$ le semigroupe engendré par $H$, c'est-à-dire l'ensemble des composées finies d'éléments de $H$. Alors deux états $x$ et $y$ sont accordables si et seulement si il existe $s \in S$ tel que $s(x) = s(y)$.

Définissons une marche aléatoire $(T_n)_{n \in \mathbf{N}}$ sur $S$ par $T_n = f_{V_0} \circ \cdots \circ f_{V_{-n+1}}$, avec la convention $T_0 = \mathrm{id}_E$. La marche aléatoire $(T_n)_{n \in \mathbf{N}}$ est reliée à la chaîne $(X_n)_{n \in \mathbf{Z}}$ par la relation $X_0 = T_n(X_{-n})$. Comme la suite des images $(T_n(E))_{n \in \mathbf{N}}$ est une suite décroissante d'ensemble finis, elle est constante à partir d'un certain rang. On notera $R_0$ la valeur limite

$$R_0 = \bigcap_{n \in \mathbf{N}} T_n(E).$$

Démontrons un premier résultat:

PROPOSITION 10. *Soit $M$ le nombre maximum d'éléments de $E$ deux-à-deux non accordables. Si $r \in S$ est un état récurrent pour la marche aléatoire $(T_n)_{n \in \mathbf{N}}$, alors $r(E)$ est formé de $M$ points deux à deux non accordables.*

PROOF. Il est clair que pour tout $s \in S$, $s(E)$ contient au moins $M$ points: les images de $M$ points de $E$ deux-à-deux non accordables. Il suffit donc de montrer que les éléments de $r(E)$ sont deux-à-deux non accordables.

Mais si $s \in S$, alors $r \circ s$ est accessible depuis $r$, donc communique avec $r$ (puisque $r$ est récurrent). Cela entraîne que

$$\operatorname{Card} r(E) = \operatorname{Card}(r \circ s)(E) \leq \operatorname{Card} s(E).$$

Par conséquent, $\operatorname{Card} r(E) = \min_{s \in S} \operatorname{Card} s(E)$. Si deux points de $r(E)$ étaient accordables par une application $t$, cela entraînerait $\operatorname{Card}(t \circ r)(E) < \operatorname{Card} r(E)$, ce qui contredirait l'égalité précédente. $\square$

REMARQUE. Comme $(T_n)_{n \in \mathbf{N}}$ est une chaîne de Markov sur un ensemble fini, elle finit presque sûrement dans une classe récurrente; La Proposition 10 montre que presque sûrement, la partie $R_0$ est formée de $M$ éléments deux-à-deux non accordables.

Venons-en au résultat principal:

THÉORÈME 11. *Presque sûrement, la loi de $X_0$ sachant $V$ est la loi uniforme sur $R_0$.*

PROOF. On utilise les résultats et une partie de la démonstration du Théorème 3. Presque sûrement, la loi de $X_n$ sachant $V_{n]} = (V_k)_{k \leq n}$ est aussi la loi de $X_n$ sachant $V$. Notons $N_n$ le nombre d'atomes de cette loi. Alors $N_n$ est aussi le nombre d'atomes de la loi de $X_{[n} = (X_k)_{k \geq n}$ sachant $V$ (puisque $X_n$ est une fonction de $X_{[n}$ et puisque $X_{[n}$ est une fonction de $X_n$ et de $V$).

Le nombre d'atomes $N_n$ est donc une fonction décroissante de $n$, mais sa loi ne dépend pas de $n$, grâce à la stationnarité du processus $(X_n, V_n)_{n \in \mathbf{Z}}$. Donc $N_n$ ne dépend pas de $N$ et est égal au nombre d'atomes $N$ de la loi de $X$ sachant $V$. D'après le Théorème 3, la loi de $X$ sachant $V$ est uniforme sur un ensemble à $N$ éléments. Comme la loi de $X_0$ sachant $V$ de déduit de la précédente et possède $N$ atomes, elle est donc aussi uniforme sur un ensemble à $N$ éléments. Comme elle est portée la partie $R_0$ (de cardinal $M$), il reste à montrer que $N \geq M$.

Pour cela, on considère (en agrandissant si besoin l'espace probabilisé) une chaîne de Markov stationnaire $(X_n^1, \ldots, X_n^M)_{n \in \mathbf{Z}}$ sur l'ensemble $A$ des $M$-uplets d'éléments de $E$ deux-à-deux non accordables, gouvernée par la relation de récurrence $X_{n+1}^k = f_{V_{n+1}}(X_n^k)$ pour tout $k \in [1, \ldots, M]$ et $n \in \mathbf{Z}$. L'existence d'une telle chaîne est assurée par le fait que $A$ est une partie absorbante finie, qui contient donc une classe récurrente. Pour tout $k \in [1, \ldots, M]$, le couple $(V, X^k)$ a même loi que le couple $(V, X)$. Comme la loi de $X$ sachant $V$ est atomique, on en déduit que presque sûrement, $X^1, \ldots, X^M$ sont des atomes de cette loi, ce qui finit la démonstration. $\square$

REMARQUE. Bien entendu, on peut remplacer l'instant 0 par tout autre instant.

3.2. *Un cas particulier intéressant.* Dans ce paragraphe, nous faisons l'*hypothèse supplémentaire* $\mathcal{H}$ suivante: il existe une probabilité $\alpha$ équivalente à $\beta$ (la loi commune des innovations $V_n$) telle que

$$\int_G \operatorname{Card} f_v^{-1}\{y\}\alpha(dv) = 1$$

pour tout $y \in E$. Ces égalités signifient que la loi uniforme sur $E$ est invariante pour les chaînes gouvernées par $f$ et par des innovations de loi $\alpha$. Sous cette hypothèse, nous allons montrer (entre autres) que le nombre maximum d'éléments de $E$ deux-à-deux non accordables divise le cardinal de $E$.

Les résultats de cette partie généralisent ceux de l'article [2], où l'on ramenait l'étude de certaines transformations du jeu de pile ou face à celle d'une chaîne de Markov constructive sur $E = \{-1,1\}^d$ gouvernée par les applications

$$f_v : (x_1, \ldots, x_d) \mapsto (x_2, \ldots, x_d, v\phi(x_1, \ldots, x_d))$$

pour $v \in \{-1,1\}$, où $\phi$ est une application fixée de $\{-1,1\}^d$ dans $\{-1,1\}$. Mais outre ce cas bien particulier, voyons un exemple simple montrant que l'hypothèse $\mathcal{H}$ est assez souvent vérifiée et un autre où elle ne l'est pas.

UN EXEMPLE CONCRET OÙ $\mathcal{H}$ EST VÉRIFIÉE. Sur un ensemble fini $E$, considérons un graphe orienté dont les arêtes sont coloriées par des couleurs numérotées de 1 à $C$ et tel que:

- De chaque point part une arête et une seule de chaque couleur.
- En chaque point arrive exactement $C$ arêtes (de couleur quelconque).

Définissons une marche aléatoire $(X_n)_{n\in\mathbf{N}}$ sur $E$ qui évolue de la façon suivante: pour passer de $X_{n-1}$ à $X_n$, on choisit d'emprunter l'arête de couleur $c \in [1, \ldots, C]$ issue de $X_{n-1}$ avec une probabilité $\beta\{c\} > 0$. La marche aléatoire ainsi définie est une chaîne de Markov constructive gouvernée par la suite des couleurs choisies et par les applications $(f_c)_{c\in[1,\ldots,C]}$ où $f_c$ est l'application dont le graphe est formé des arêtes de couleur $c$. Dans ce cas, l'hypothèse $\mathcal{H}$ est vérifiée en prenant $\alpha$ égale à la loi uniforme sur $[1, \ldots, C]$.
□

UN EXEMPLE OÙ $\mathcal{H}$ N'EST PAS VÉRIFIÉE. On prend $E = \{1,2,3,4\}$, $G = \{1,2,3\}$ et $H = \{f_1, f_2, f_3\}$ avec

$$f_1(1) = 2, \quad f_1(2) = 3 \quad \text{et} \quad f_1(3) = f_1(4) = 1;$$
$$f_2(1) = 2, \quad f_2(2) = 4 \quad \text{et} \quad f_2(3) = f_2(4) = 1;$$
$$f_3(1) = 1, \quad f_3(2) = 2 \quad \text{et} \quad f_3(3) = f_3(4) = 3.$$

On vérifie facilement que la chaîne de Markov constructive associée aux applications $f_1$, $f_2$ et $f_3$ est irréductible et apériodique. Comme $f_v^{-1}\{4\}$ a un cardinal toujours strictement inférieur à $f_v^{-1}\{1\}$, l'hypothèse $\mathcal{H}$ n'est pas vérifiée. D'ailleurs, le semigroupe engendré par $H$ est $S = \{f_1, f_2, f_3, f_1^2, f_2^2\}$ et les seuls états accordables sont 3 et 4. Le nombre maximum d'états deux-à-deux non accordables est 3, et ne divise donc pas le cardinal de $E$.

LES RÉSULTATS OBTENUS GRÂCE À L'HYPOTHÈSE $\mathcal{H}$. Notons $S$ le semi-groupe engendré par $H = \{g \in E^E : P[f_{V_1} = h] > 0\}$. Nous dirons que des états $a_1, \ldots, a_m$ sont simultanément accordables s'il existe $s \in S$ tel que $s(a_1) = \cdots = s(a_m)$. Notons $N$ le nombre maximal d'éléments simultanément accordables. Énonçons un premier résultat:

LEMME 12. *Soit $F$ une partie de $E$ formée de $N$ éléments simultanément accordables. Alors toute image réciproque de $F$ par une application de $S$ est formée de $N$ éléments simultanément accordables.*

PROOF. Il est clair que l'image réciproque de $F$ par une application de $S$ est formée d'états simultanément accordables, donc d'au plus $N$ états. Mais en sommant sur $y \in F$ les égalités $\int_G \operatorname{Card} f_v^{-1}\{y\} \alpha(dv) = 1$ on obtient $\int_G \operatorname{Card} f_v^{-1}(F) \alpha(dv) = N$. Mais comme $\operatorname{Card}(f_v^{-1}(F)) \leq N$, on a donc $\operatorname{Card}(f_v^{-1}(F)) = N$ pour $\beta$-presque tout $v \in G$. Autrement dit, $\operatorname{Card}(h^{-1}(F)) = N$ pour $h \in H$. Le résultat s'étend immédiatement à toute composée d'éléments de $H$.

Le Lemme 12 permet de montrer un résultat remarquable sur les composées réalisant un accord maximal, c'est-à-dire les applications $s \in S$ telles que $s(E)$ soit formée de $M$ points deux-à-deux non accordables.

THÉORÈME 13. *Si $s \in S$, réalise un accord maximal, alors les images réciproques par $s$ des éléments de $s(E)$ contiennent chacune $N$ éléments (simultanément accordables). On a donc $MN = \operatorname{Card} E$.*

*En particulier, si $\operatorname{Card} E$ est premier, alors:*

- *Soit tout élément de $H$ est une bijection de $E$ dans $E$ et pour tout $n \in \mathbf{Z}$, $X_n$ est indépendante de la suite $V$ et de loi uniforme sur $E$.*
- *Soit les éléments de $E$ sont simultanément accordables, et $\mathcal{F}_n^X \subset \mathcal{F}_n^V$ pour tout $n \in \mathbf{Z}$.*

PROOF. Comme les images réciproques par $s$ des éléments de $s(E)$ sont formées chacune d'au plus $N$ éléments simultanément accordables, il suffit de montrer l'égalité $MN = \operatorname{Card} E$ pour conclure. Il suffit donc de construire une application pour laquelle le résultat du théorème est vérifié. Cette construction découle de l'application répétée du lemme ci-dessous. □

LEMME 14. *Soient $F_1, \ldots, F_k$ des parties disjointes de $E_d$ formées chacune de $N$ états simultanément accordables, et $s \in S$ une application constante sur chacune de ces parties. Notons $c_1, \ldots, c_k$ la valeur prise par $s$ sur les parties $F_1, \ldots, F_k$. Alors:*

- *Les éléments $c_1, \ldots, c_k$ sont deux-à-deux non accordables, d'où $k \leq M$.*
- *Si $F_1 \cup \cdots \cup F_k \neq E$ (en particulier, si $k < M$), on peut construire $k+1$ parties disjointes de $E_d$ formée chacune de $N$ éléments simultanément accordables, et une application de $S$ constante sur chacune de ces parties.*

PROOF. Le premier point est immédiat: s'il existait deux indices $i < j$ tels que les points $c_i$ et $c_j$ soient accordables par une composée $t$, alors la composée $t \circ s$ accorderait les $2M$ éléments de $F_i \cup F_j$. Montrons le second point: soit donc $a \in E_d \setminus (F_1 \cup \cdots \cup F_k)$ et $c_{k+1} = s(a)$. Alors $c_{k+1}$ est différent de $c_1, \ldots, c_k$ (sinon, on obtiendrait $N+1$ éléments simultanément accordables en ajoutant $a$ à l'une des parties $F_1, \ldots, F_k$). Soit $t \in S$ une application envoyant $c_1$ sur $a$. Les images réciproques de $c_1, \ldots, c_k, c_{k+1}$ par $s \circ t \circ s$ sont formées chacune de $N$ éléments simultanément accordables. □

**4. Addendum: une généralisation d'Aldous et Bandyopadhyay.** Depuis la soumission de cet article, Aldous et Bandyopadhyay ont publié dans [1] un résultat tout à fait similaire et plus général que notre Théorème 9. Ces auteurs se placent dans le cadre de processus récursifs indexés par les sommets d'un arbre de Galton–Watson, qui généralisent les chaînes de Markov constructives indexées par $-\mathbf{N}$.

Ils obtiennent dans ce cadre une équivalence (endogeny $\Leftrightarrow$ bivariate uniqueness property) semblable à l'équivalence (2) $\Leftrightarrow$ (3) de notre Théorème 9, ainsi que l'équivalence avec une notion qui peut s'apparenter avec l'accordabilité, mais sur un espace d'états quelconque (polonais).

Il convient de noter, toutefois, que leur démonstration de l'implication bivariate uniqueness property $\Rightarrow$ endogeny nécessite une hypothèse de continuité dont on peut se passer par une adaptation immédiate de notre démonstration de l'implication (3) $\Rightarrow$ (2), ce qui répond à la question ouverte 12 de [1].

LABORATOIRE DE MATHMATIQUES
INSTITUT FOURIER
UMR5582 (UJF-CNRS)
BP 74
38402 ST MARTIN D'HÈRES CEDEX
FRANCE
E-MAIL: Jean.Brossard@ujf-grenoble.fr
         Christophe.Leuridan@ujf-grenoble.fr